\documentstyle[sprocl]{article}
\input{psfig}
\def\be{\begin{equation}}
\def\ee{\end{equation}}
\def\bea{\begin{eqnarray}}
\def\eea{\end{eqnarray}}

\begin{document}
\title{ON THE FUNDUMENTAL INVARIANT OF THE HECKE ALGEBRA $H_{n}(q)$}
\author{J. Katriel$^{\dagger}$, B. Abdesselam$^{\ddagger}$ and 
A. Chakrabarti$^{\ddagger,}$
\footnote{ \it Talk presented at the Nankai workshop, Tianjin, 1995 by 
A. Chakrabarti.}}
\address{$^{\dagger}$ Department of Chemistry, Technion, 32000 Haifa, 
Israel.} 
\address{$^{\ddagger}$Centre de Physique Th{\'e}orique, 
Ecole Polytechnique, 
91128 Palaiseau Cedex, France.
Laboratoire Propre du CNRS UPR A.0014}


\maketitle
\abstracts{ The fundumental invariant of the Hecke algebra $H_{n}(q)$ is the 
$q$-deformed class-sum of transpositions of the symmetric group $S_{n}$. 
Irreducible representations of $H_{n}(q)$, for generic $q$, are shown to be 
completely characterized by the corresponding eigenvalues of $C_{n}$ alone. 
For $S_{n}$ more and more invariants are necessary as $n$ inereases. It 
is pointed out that the $q$-deformed classical quadratic Casimir of $SU(N)$ 
plays an analogous role. It is indicated why and how this should be a general 
phenomenon associated with $q$-deformation of classical algebras. Apart from 
this remarkable conceptual aspect $C_{n}$ can provide powerful and elegant 
techniques for computations. This is illustrated by using the sequence 
$C_{2}$, $C_{3}, \cdots,\; C_{n}$ to compute the characters of $H_{n}(q)$.}

This talk will be based on $[1]$ and $[2]$. Much more complete discussions 
and references to other authors can be found there. More recent developments 
can be found in $[3]$. 

Let me start by recapitulating certain facts concerning the invariants of the 
classical symmetric group $S_{n}$. The single cycle class-sums 
$\lbrace [p]_{n};\;p=2,\;3,\cdots, \;n \rbrace$ belong to the centre. Here 
$[2]_{n}$ is the sum of transpositions, $[3]_{n}$  is that of circular 
permutation of triplets (each term being a product of transpositions) and so 
on. Their eigenvalues characterize irreducible representations (irreps.) of 
$S_{n}$ corresponding to different standard Young tableaux with $n$ boxes.

\proclaim Definition. {\it Content} of the box in the $i$-th row and $j$-th 
column of the Young tableau is equal to $(j-i)$.   

The symmetric power sums of these contents gives the eigenvalues 
$\lambda_{[p]_{n}}^{\Gamma}$ of $[p]_{n}$ for $\Gamma$. Thus 
\begin{equation}
\begin{array}{l}
\lambda_{[2]_{n}}^{\Gamma}=\displaystyle\sum_{(i,\;j) \in \Gamma} (j-i)
\end{array}
\end{equation}
\begin{equation}
\begin{array}{l}
\lambda_{[3]_{n}}^{\Gamma}=\displaystyle\sum_{(i,\;j) \in \Gamma} (j-i)^{2}-
{1 \over 2}n(n-1)
\end{array}
\end{equation}
and so on $[1]$.

For $n \geq 6$ the eigenvalues of $[2]_{n}$ show degeneracy. As $n$ inereases 
higher and higher $[p]_{n}$'s are needed to uniquely characterize each 
irreducible representation.

{\it The situation changes dramatically as $S_{n}$ is $q$-deformed to 
$H_{n}(q)$. The eigenvalues of $C_{n}$, the $q$-deformed $[2]_{n}$, alone 
suffice to characterize the irreducible representations for arbitrary $n$. 
(Throughout only real, positive i.e. generic $q$ is considered.)}
 
I will now show how this becomes possible. The generators of the $H_{n}(q)$ 
satisfy 
\begin{equation}
\label{equation:Hecke}
\begin{array}{ll}
g_i^2=(q-1)g_i+q &\;\;\; i=1,\, 2,\, \cdots,\, n-1 \\
g_ig_{i+1}g_i=g_{i+1}g_ig_{i+1} &\;\;\; i=1,\, 2,\, \cdots,\, n-2 \\
g_ig_j=g_jg_i &\;\;\; {\mbox{if }} |i-j|\geq 2
\end{array}
\end{equation}
For $q=1$ one gets $S_{n}$. The fundumental invariant is
\begin{eqnarray}
\label{equation:fund}
C_n &=& g_1+g_2+\cdots +g_{n-1}+
\frac{1}{q}(g_1g_2g_1+g_2g_3g_2+\cdots +g_{n-2}g_{n-1}g_{n-2}) \nonumber \\ &+&
\frac{1}{q^2}(g_1g_2g_3g_2g_1+g_2g_3g_4g_3g_2+\cdots +g_{n-3}g_{n-2}g_{n-1}g_{n-2}g_{n-3})
\nonumber \\ &+&\cdots \nonumber \\ &+&
\frac{1}{q^{n-2}}g_1g_2\cdots g_{n-2}g_{n-1}g_{n-2}\cdots g_2g_1 
\end{eqnarray}
For $q=1$, $g_{1}g_{2}g_{1}=(13)$ and so on and one gets back $[2]_{n}$.
The eigenvalue of the fundumental invariant for the $Y$-tableau $\Gamma$ can be 
shown $[1]$ to be the following $q$-deformation of $(1)$,
\begin{equation}
\label{equation:eigen}
\Lambda_n^{\Gamma}=q\displaystyle\sum_{(i,j)\in\Gamma}{q^{j-i}-1 \over q-1}
=q\displaystyle\sum_{(i,j)\in\Gamma} [j-i]_q.
\end{equation}

\proclaim Definition. $q$-content of the box in the $i$-th row and $j$-th 
column of the Young tableau is equal to $q\;[j-i]_{q}$.

Hence $\Lambda_n^{\Gamma}$ is the sum of the $q$-contents of the boxes of 
$\Gamma$. 

Consider, for $n=6$, the irreducible representations $[4,\;1,\;1]$ and 
$[3,\;3]$. The box contents are

\medskip
\centerline{\psfig{figure=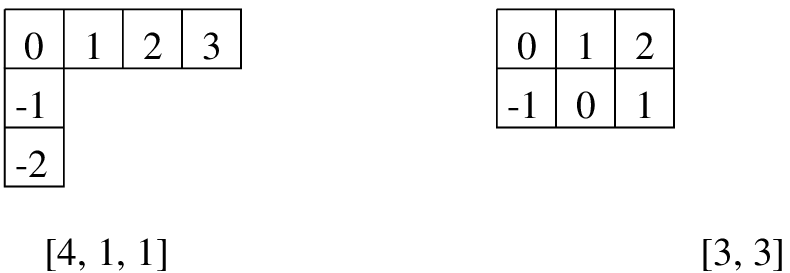}}
\medskip

For $S_{6}$,
\begin{equation}
\Lambda_6^{[4,\;1,\;1]}=\Lambda_6^{[3,\;3]}=3
\end{equation}

For $H_{6}(q)$,
\begin{equation}
\begin{array}{l}
\Lambda_6^{[4,\;1,\;1]}(q)=q^3+2\;q^2+3\;q-2-{1 \over q} \\
\Lambda_6^{[3,\;3]}(q)=q^2+3 q -1
\end{array}
\end{equation}
{\it Thus the degeneracy is lifted as $q$ moves away from unity.} This is the 
simplest non-trivial example. For the general case one notes:

(i) The $q$-contents are constant for boxes on the same diagonal of
$\Gamma$.

(ii) Developping the $q$-brackets and regrouping terms 
\begin{equation}
\Lambda_{n}^{\Gamma}=\displaystyle\sum_{k >0 }q^{k} \pi_{k}^{\Gamma}-
\displaystyle\sum_{k < 0 }q^{k+1} \nu_{k}^{\Gamma}
\end{equation}
where
$$
\pi_{k}^{\Gamma}=\displaystyle\sum_{l \geq k >0 }(\hbox{number of boxes with content $l$})
$$
$$
\nu_{k}^{\Gamma}=\displaystyle\sum_{l \leq k <0 }(\hbox{number of boxes with content $l$})
$$
From $(i)$ and $(ii)$ it is not difficult to show $[1]$ that 
{\it $\Lambda_{n}^{\Gamma}$ completely determines $\Gamma$ and hence the 
irrep.}

Set $q=e^{\delta}\;(\not = 1)$ and let
\begin{equation}
\tilde{C}_{n}\equiv {q-1 \over q} C_{n}
\end{equation}
then $[1]$,
\begin{equation}
\Lambda_{\tilde{C}_{n}}^{\Gamma}=\delta \lambda_{[2]_{n}}^{\Gamma}+
{\delta^2 \over 2} ( \lambda_{[3]_{n}}^{\Gamma}+ {1\over 2}n(n-1))+\cdots 
\end{equation}
The eigenvalues of all $[p]_{n}$ ($p=2,\cdots,\;n$) in the coefficients
of the above series. This is another way of exhibiting that $C_{n}$ by itself 
contains information equivalent to that supplied by all $[p]_{n}$ for $S_{n}$. 

Projection operators for irreps. can be constructed in terms of $C_{n}$ in a 
straightforward way since there is no degeneracy. When the limit 
$q\rightarrow 1$ is taken correctly higher class-sums of $S_{n}$ appear 
automatically as necessary to project out the corresponding irrep. of $S_{n}$. 
This is discussed in detail in $[1]$. Further interesting uses of projection 
operators can be found in $[3]$.

In $[1]$ a direct relation was given between the eigenvalues of $C_{n}$ an 
those of the Casimir of $SU_{q}(N)$ ($q$-deformation of the Casimir quadratic 
in the Cartan-Weyl generators of $SU(N)$) for an irrep. corresponding to a 
$Y$-diagram $\Gamma$ with $n$ boxes (and at most $N-1$ rows). It was shown 
that this Casimir $C_{2}$ can be so redefined (denoted then by 
$\tilde{C}_{2}$) that the eigenvalue for $\Gamma$ is just
\begin{equation}
\Lambda_{\tilde{C}_{2}}^{\Gamma}=\displaystyle\sum_{k=1}^{N-1}q^{2(l_{k}-k)}
\end{equation}
where $l_{k}$ is the number of boxes in the $k$-th row. This was derived 
using the Gelfand-Zetlin basis $[1]$. But $(11)$ is, of course, independent 
of the choice of such a basis. Since
\begin{equation}
l_{k} \geq l_{k+1},\;\;\;\;\;\;\;\;\;(l_{k}-k) > (l_{k+1}-k-1)
\end{equation}
Hence the indices of $q$ in $(11)$ are strictly monotonically decreasing.

Thus even for a reducible representation arising in a certain context 
(say some model) if one obtains the matrix of $\tilde{C}_{2}$ from some source 
and diagonalizes it the coefficient of each block of unit matrix must be of 
the form
\begin{equation}
\displaystyle\sum_{k=1}^{N-1}q^{2 L_{k}} \;\;\;\;\;\;\;(L_{k} > L_{k+1})
\end{equation}
Now setting 
\begin{equation}
l_{k}=L_{k}+k,\;\;\;\;\;\;\;\;\;(k=1,\cdots, \;N-1)
\end{equation}

The $Y$-diagram $\Gamma$ is completely determined. {\it Thus the eigenvalue 
of a suitably $q$-deformed quadratic Casimir completely characterizes an 
irrep. of $SU_{q}(N)$.} (For $q=1$ or $SU(N)$ one needs, in general, all 
the invariants upto order $N$.)

Setting $q=e^{\delta}$,
\begin{equation}
\displaystyle\sum_{k}q^{2 L_{k}}=1+(2\delta)(\displaystyle\sum_{k}L_{k})+
{1\over 2!}(2\delta)^{2}(\displaystyle\sum_{k}L_{k}^{2})+\cdots
\end{equation}
One can compare $(15)$ with $(10)$. 

The coefficients of higher powers of $\delta$ contain informations equivalent 
to those of higher order Casimirs of $SU(N)$. 

I present now, without derivation, the relation between $C_{n}$ and 
$\tilde{C}_{2}$ eigenvalues for a $\Gamma$ with $n$ boxes $[1]$,
\begin{equation}
\left({q^{2}-1 \over q^2}\right)^{2}\Lambda_{C_{n}(q^2)}^{\Gamma}+
{q^{2}-1 \over q^2}n =\Lambda_{\tilde{C}_{2}}^{\Gamma}+{q^{2(-N+1)}-1 \over q^2 -1}
\end{equation}
for ($\displaystyle\sum l_{k}=n$ in $(11)$).

The Hecke and $q$-deformed unitary algebras are well-known to be closely 
related. But the aspect presented here is more general. {\it Thus the 
$q$-deformed quadratic Casimirs of the other Lie algebras should play 
analogous roles}. Our investigation is not complete. Here only $SU_{q}(N)$ has 
been studied. However in an accompanying talk $[4]$ the foregoing statement 
is confirmed for $SO_{q}(5)$. The discussion at the end of $[4]$ gives an 
idea of the richness of content of the $q$-deformed Casimirs. 

Apart from such remarkable conceptual aspects, $C_{n}$ or, even better, the 
sequence ($C_{2}$, $C_{3},\cdots,\;C_{n}$) nested in $H_{n}(q)$ can furnish 
powerful techniques for various goals. As an example, I will indicate below 
how they can be used to compute characters. A detailed study can be found in 
$[2]$. (Another interesting aspect has been studied in $[3]$). 

For the sequence $H_{2}(q)\subset H_{3}(q) \subset \cdots \subset H_{n}(q)$ 
one defines the Murphy operators
\begin{equation}
L_{2}=C_{2}, \;\;L_{3}=C_{3}-C_{2},\cdots , \;\;L_{n}=C_{n}-C_{n-1}
\end{equation}
One obtains
\begin{equation} 
L_{p}=\displaystyle\sum_{i=1}^{p-1} q^{1-p+i}(g_{i} g_{i+1}\cdots
g_{p-1}\cdots g_{i+1} g_{i}) 
\end{equation}
and
\begin{equation} 
L_{p+1}={1 \over q} g_{p}L_{p}g_{p} +g_{p}
\end{equation}

Basis vectors of an irrep. can be specified by sequences of 
$Y$-diagrams (indicating how successive boxes are added)
\begin{equation} 
\Gamma_{2} \subset \Gamma_{3} \subset \cdots \subset \Gamma_{n}
\end{equation}
The eigenvalue of $L_{i}$ can be shown to be $[2]$
\begin{equation} 
\lbrace \Gamma_{i}\backslash\Gamma_{i-1}\rbrace_{q}\equiv q[k_{i}-p_{i}]_{q}
\end{equation}
where $\Gamma_{i}$ is obtained by adding the box $(k_{i},\;p_{i}$)
to $\Gamma_{i-1}$. The eigenvalue is the $q$-content of the last 
box added. Also
\begin{equation} 
tr(L_i)_{\Gamma_n}=\displaystyle\sum_{\Gamma_{n-1}\subset\Gamma_n} 
tr(L_i)_{\Gamma_{n-1}}\;\;\; (i=2,\, 3,\, \cdots,\, n-1.)
\end{equation}
and
\begin{equation} 
tr(L_n)_{\Gamma_n}= 
\displaystyle\sum_{\Gamma_{n-1}\subset\Gamma_n} |\Gamma_{n-1}|\, 
\{ \Gamma_n\setminus\Gamma_{n-1} \}_q\; .
\end{equation}
where
\begin{equation} 
|\Gamma_{n-1}|=\displaystyle\sum_{\Gamma_{n-2}\subset\Gamma_{n-1}} 
|\Gamma_{n-2}|=dim\;\Gamma_{n-1}
\end{equation}
For what follows we will need traces of products of 
{\it non-consecutive} Murphy operators only. For such products 
with
$$
\alpha_{i+1} \geq \alpha_{i}+2
$$
\begin{equation} 
tr\left(\prod_{i=1}^{\ell} L_{\alpha_i}\right)_{\Gamma_n} =
\displaystyle\sum_{\Gamma_{n-1}\subset\Gamma_n} 
tr\left(\prod_{i=1}^{\ell} L_{\alpha_i}\right)_{\Gamma_{n-1}}\;\;\;(\hbox{for}\;\alpha_{l} < n) 
\end{equation}
\begin{equation} 
tr\left(\prod_{i=1}^{\ell} L_{\alpha_i}\right)_{\Gamma_n} =
\displaystyle\sum_{\Gamma_{n-1}\subset\Gamma_n} 
\{\Gamma_n\setminus\Gamma_{n-1}\}_q \;
tr\left(\prod_{i=1}^{\ell-1} L_{\alpha_i}\right)_{\Gamma_{n-1}}\;\;(\hbox{for}\; \alpha_{l}=n)
\end{equation}
The recursion relations $(22)$ to $(26)$ yield easily the traces 
of the $L$'s and their non-consecutive products $[2]$. A symbolic
program is easy to set up. Taking the trace of each side of $(18)$ and 
inverting the relation one obtains $[2]$. 
\begin{equation}
tr\Big(g_1 g_2 \cdots g_{k-1} \Big)=\left(\frac{q}{q-1}
\right)^{k-2}\displaystyle\sum_{i=0}^{k-2} (-1)^i{{k-1}\choose i} tr(L_{k-i})
\end{equation}
Similarly, after multiplying both sides of $(18)$ by $L_{m}$ (non-consecutive),
\begin{equation}
tr\Big((g_1 g_2 \cdots g_{k-1}) L_m\Big)=\left(\frac{q}{q-1}\right)^{k-2} \displaystyle\sum_{i=0}^{k-2} (-1)^i 
{{k-1}\choose{i}}
tr(L_{k-i} \; L_m) \; .
\end{equation}
Continuing step-wise, with suitable choices of $m$ at each step, it can be 
shown $[2]$ that one finally obtains in terms of known traces of the type 
$(25)$ and $(26)$ traces of the form
\begin{equation}
tr\Big((g_1 g_2 \cdots g_{m_{1}-1})(g_{m_{1}+1} \cdots 
g_{m_{2}-1}) \cdots (g_{m_{j}+1} \cdots 
g_{p})(g_{k}g_{k+1}\cdots g_{p+r}\cdots g_{k+1}g_{k})\Big) 
\end{equation}
Here the last factor comes from the last non-consecutive $L$ ($L_{p+r-1}$). 
This, in general, has overlapping indices with the preceding factors. At 
each previous step such overlaps has been assumed to be {\it reduced} 
(reexpressed as sums of traces of ordered products the $g_{i}$'s in ascending 
order of $i$) so that they have no overlap but, possibly, {\it cuts} (at $i=m_{1}, \;m_{2},\cdots, \;m_{j}, \hbox{say})$. This reduction procedure, to be 
applied again to $(29)$, will be briefly described below. Let us however first 
are exhibit the simplest results, illustrating general properties. One obtains $[2]$
\begin{equation}
tr(g_{i})=tr(L_{2}),\;\;\;\;\;\;\;(i=1,\cdots,\;n-1) 
\end{equation}
\begin{equation}
tr(g_{i}g_{i+1})=
\Big({q \over q-1}\Big)\Big(tr(L_{3})-2tr(L_{2})\Big),\;\;\;\;\;
(i=1,\cdots,\;n-2)  
\end{equation}
\begin{equation}
tr(g_{i}g_{i+2})={1\over q-1}\Big( -2q\, tr(L_{2})+
(q+1)^{2}tr(L_{3})-(1+q^{2})tr(L_{4})+(q-1)tr(L_{2}L_{4})\Big) 
\end{equation}
The equalities of the traces in each example illustrate a fundumental lemma 
$[2]$: {\it The trace of product of any number of disjoint
sequences, in any irrep., depends only on the lengths of the component 
connected sequences.}

Note that in $(32)$, $(g_{i}g_{i+2}$) being disjoint (i.e. a cut at $i+1$) a 
non-consecutive product $L_{2}L_{4}$ appears on the right. {\it Note that 
the $L$'s on the right do not depend on $n$ (of $H_{n}(q)$)}. This is also 
a general feature.

When there is no overlap and at most one cut $(29)$ can be reduced relatively 
easily $[2]$. Thus
\begin{equation}
\begin{array}{ll}
V_k &\equiv tr\Big((g_{1}...g_{k-1})(g_{k+1}...g_{p})(g_{p+1}...
g_{p+r}...g_{p+1})\Big) \\
&= \displaystyle\sum_{\ell=0}^{r-1}{{r-1}\choose l}q^{l}(q-1)^{r-l-1}
tr\Big((g_{1}...g_{k-1})(g_{k+1}...g_{p+r-l})\Big)
\end{array}
\end{equation} 
(For $V_{0}$ with $k=0$ the first factor is defined to be unity). When 
there is an overlap we introduce {\it $f$-expansions} defined below. (For 
a full account see Appendix $[2]$). From,
$$
g_{i}^{2}=(q-1)g_{i}+q
$$
one deduces 
\begin{equation}
g_{i}^{p}=f_{p}g_{i}+q\;f_{p-1}
\end{equation} 
where
$$
f_{p}={q^{p}-(-1)^{p} \over q+1 } 
$$
It can be shown that (for overlap$=p-l+1$)
\begin{equation}
\label{equation:e14}
tr\Big((g_{1}...g_{p})(g_{l}...g_{p+r}...g_{l})\Big)
=(q-1)\displaystyle\sum_{r=1}^{p-l+1} q^{k}f_{2(p-l+1-k)+1}V_{k}
+f_{2(p-l+1)+1}V_{0}
\end{equation} 

{\it Thus the $f$-coefficients are determined only by the length of the 
overlap.} The general case with overlap and multiple cuts is treated in 
App. $[2]$. 

Tables of characters (polynomials in $q$) are given in $[2]$. Here let me 
just summerize the main steps:

(1) Traces of (non-consecutive) products of Murphy operators.

(2) Traces of products of $g$'s with cuts and overlap in terms of $(1)$. 

(3) Reduction removing overlaps ($f$-expansions).

\end{document}